\documentclass[matzei,envcountsame,draft]{svjour}

\usepackage{geometry}                % See geometry.pdf to learn the layout options. There are lots.
\usepackage{graphicx}
\usepackage{amsmath,amssymb,latexsym}
\usepackage{epstopdf}
	\def\copyleft{}
	\def\makeheadbox{{%
	\hbox to0pt{\vbox{\baselineskip=10dd\hrule\hbox
	to\hsize{\vrule\kern3pt\vbox{\kern3pt
	\hbox{
	submitted version}
	\kern3pt}\hfil\kern3pt\vrule}\hrule}%
	\hss}}}

\def\Romannumeral#1 {\begingroup{\uppercase\expandafter{\romannumeral#1}}\endgroup}% 
\def\pf#1{{\def\temp{#1} 
              \ifx\temp\empty 
                  \noindent\slshape\textbf{Proof.\ }
              \else 
                  \noindent\slshape\textbf{Proof of\ #1.\ }
              \fi}}
\def\qed{\hspace*{\fill}$\square$}
\def\Stab#1_#2.{#2_{\{#1\}}}%
\def\Fix#1_#2.{{#2_{#1}}}%
\def\Aut#1_#2.{\mathsf{Aut}_{#2}(#1)}%
\def\Isom#1_#2.{\mathsf{Isom}_{#2}(#1)}%
\def\bs#1.{
              \def\temp{#1} 
              \ifx\temp\empty 
                   \mathcal{B}
              \else
                   \mathcal{B}(#1)
              \fi
}

\renewcommand{\emph}{\textbf}%

\newcommand{\svol}{\mathop\mathrm{vol}}
\newcommand{\Svol}{\mathop{s\text{-}\mathrm{vol}}}

\spnewtheorem{algorithm}[theorem]{Algorithm}{\bf}{\sf}

\title{Scales for co-compact embeddings of virtually free groups}
\author{%workgroup tdlcG} % 
Udo Baumgartner} 
\institute{
School of Mathematical and Physical Sciences, 
The University of Newcastle, University Drive, Building V, 
Callaghan, NSW 2308, Australia\\
\email {Udo.Baumgartner@newcastle.edu.au} %, George.Willis@newcastle.edu.au}
}
\date{\today}                                           % Activate to display a given date or no date

\begin{document}
\maketitle
\begin{abstract}
Let $\Gamma$ be 
a group 
which is 
virtually free 
of rank at least $2$ 
and 
let $\mathcal{F}_{td}(\Gamma)$ be 
the family 
of totally disconnected, locally compact groups 
containing $\Gamma$ 
as a co-compact lattice.  
We prove that 
the values of 
the scale function 
with respect to groups in $\mathcal{F}_{td}(\Gamma)$ 
evaluated on the subset $\Gamma$ 
have only finitely many 
prime divisors. 
This %result 
can be thought of as 
a uniform property of 
the family $\mathcal{F}_{td}(\Gamma)$. 
 


\end{abstract}

\section{Introduction}

%\nocite{Mostow-Margulis-rig-lc-targets}

Let $\Gamma$ be 
a finitely generated group 
and let $S$ be 
a finite set of generators 
for $\Gamma$. 
The group $\Gamma$ 
embeds 
as a cocompact, 
discrete subgroup 
in the automorphism group 
of its Cayley graph 
with respect to $S$, 
which is 
a totally disconnected, locally compact group. 

This note 
begins the examination of  
common features of 
the class $\mathcal{F}(\Gamma)$ 
of all embeddings  
of a fixed finitely generated group $\Gamma$ 
as a discrete, cocompact subgroup 
in some locally compact group. 
By the previous paragraph,    
the automorphism group 
of any 
Cayley graph  
of $\Gamma$ 
belongs to $\mathcal{F}(\Gamma)$. 
In this article 
we treat 
the case 
where  
$\Gamma$ 
is virtually free 
of finite rank 
at least $2$. 
Most of our results 
will be about 
the subclass $\mathcal{F}_{td}(\Gamma)$ 
of $\mathcal{F}(\Gamma)$ 
consisting of 
cocompact embeddings 
of $\Gamma$ 
into totally disconnected, locally compact groups; 
see however Corollary~\ref{cor:tdlcG-envelope(freeG)}.

This program 
was suggested 
(in a less general form)  
by George Willis 
in~ \cite[penultimate topic in Section~6]{can_form(aut(tdlcGs))}.  
One of the common features 
of the class $\mathcal{F}_{td}(\Gamma)$ 
suggested there 
for examination 
is the set 
of values of the scale function,  
whose definition we recall below, 
with respect to 
%totally disconnected 
members of $\mathcal{F}_{td}(\Gamma)$ 
evaluated on $\Gamma$. 

%
%Alexander 
Furman 
undertook 
a related project 
for lattices in semisimple (connected real) Lie groups 
in~\cite{Mostow-Margulis-rig-lc-targets}. 
He 
proposed 
to classify 
all second countable, locally compact groups,  
which admit 
a lattice embedding
(not necessarily cocompact) 
of a given lattice $\Gamma$ 
in a semisimple connected, real Lie group 
(%see 
%Problem 
%on page~31 in~\cite{Mostow-Margulis-rig-lc-targets}) 
\cite[p.~31]{Mostow-Margulis-rig-lc-targets})
and 
solved this problem 
if $\Gamma$ is 
a lattice  
in a simple Lie group  
of higher rank 
(\cite{Mostow-Margulis-rig-lc-targets}, Theorem~A). 
He was also able 
to classify 
second countable, locally compact groups 
which admit 
lattice embeddings  
which are 
cocompact,  
if $\Gamma$ is 
an irreducible lattice  
in either 
a semisimple, connected, real Lie group  
not locally isomorphic to $\mathrm{SL}_2(\mathbb{R})$ 
or 
a cocompact lattice 
in a group 
locally isomorphic to  $\mathrm{SL}_2(\mathbb{R})$ 
(\cite{Mostow-Margulis-rig-lc-targets}, 
Theorem~B and Theorem~C).

We now 
return to 
the question, 
how one might 
characterize 
common features 
of discrete, cocompact embeddings 
of a given group $\Gamma$ 
into some 
totally disconnected, locally compact group. 
One way 
in which 
this may be done, 
following the suggestion 
by Willis 
mentioned above,  
is by restricting 
the values 
taken by 
the scale functions 
with respect to 
the codomain 
of the embeddings 
on the image of $\Gamma$.  
 
For 
a totally disconnected, locally compact group $G$, 
the value of 
the \emph{scale function $s_G$} 
at an automorphism, 
$\alpha$,  
of $G$ 
measures 
minimal distortion 
of compact, open subgroups 
of $G$ 
under $\alpha$. 
The scale function 
is defined on 
the set of automorphisms of $G$ by 
the formula 
\[
s_G(\alpha):=\min\{|\alpha(V)\colon\alpha(V)\cap V|\colon V\text{ a compact, open subgroup of }G\}\,.
\]
Note that 
the minimum 
above 
is attained, 
because 
it is 
%applied to 
formed for 
a set of 
positive integers. 
A compact, open subgroup $O$ 
of $G$ 
is \emph{tidy for $\boldsymbol\alpha$} 
if this minimum is attained at $O$. 
The scale of 
a group element, 
$x$, 
of $G$ 
is defined to be 
the value of 
the scale function 
with respect to $G$ 
at conjugation by $x$. 

The collection of tidy subgroups 
for an automorphism 
and the scale function 
are invariants, 
which have been used 
to answer various questions 
on totally disconnected, locally compact groups; 
see~\cite{can_form(aut(tdlcGs))} 
for a survey. 

The scale of an automorphism %$\alpha$ 
of a locally compact, totally disconnected group %$G$  
is an analogue 
for the set of eigenvalues 
of a linear transformation. 
%Indeed, 
%when $G$ is a Lie group, 
%the image of an automorphism $\alpha$ 
%under the adjoint representation 
%in case $\alpha$ is 
%an automorphism 
%of a Lie group.  
%
That there should be 
any uniform bound 
on the primes 
dividing 
values 
of the scale 
with respect to 
the elements of $\mathcal{F}_{td}(\Gamma)$ 
is not clear, 
even though 
a single element of $\mathcal{F}_{td}(\Gamma)$  
contributes only 
a finite number 
of prime factors 
by Theorem~3.4 
in~\cite{prime-factors(sF(cp.gen))=finite}. 

The author 
was surprised 
to discover, 
that 
even a slightly stronger result 
can be proved 
quite easily 
in the case when 
$\Gamma$ is virtually 
a free group 
of rank at least $2$, 
by putting together 
work by 
Lee Mosher, Michah Sageev, and Kevin Whyte 
and Alexander Lubotzky; 
see Corollary~\ref{cor:bound(div(scales))_virt-freeGs}.  

It can be shown,  
that 
virtually abelian groups 
embed cocompactly 
only 
in totally disconnected, locally compact groups 
whose scale function 
is identically $1$; 
hence 
a result 
analogous to Corollary~\ref{cor:bound(div(scales))_virt-freeGs}  
does hold 
for virtually infinite cyclic groups 
also. 
As a further remark, 
by the results 
of Furman  
mentioned above,  
groups $\Gamma$ 
to which 
one of his Theorems~A, B, or~C 
applies,  
also 
embed cocompactly only 
in a totally disconnected, locally compact group $G$ 
whose scale function 
is identically~$1$, 
\underline{provided}  
the group~$G$ 
is second countable.

The method of proof used 
in this note 
suggests, 
that 
replacing 
the use of Theorem~\ref{thm:virtF,cocp latt} 
by 
an appeal to the main result of~\cite{large_scale-Geo(prod(trees))} 
will prove 
a generalization 
of Corollary~\ref{cor:bound(div(scales))_virt-freeGs} 
to groups $\Gamma$, 
which are quasi-isometric to  
a product of finitely many trees 
once 
a replacement for Lemma~\ref{lem:ramification_bound(SchottkyGs)} 
can be found,  
to bound 
the quotient 
of an action 
of a group $\Gamma$ 
on a product of trees 
in terms of properties of $\Gamma$. 

Note that  
the class of groups 
which are quasi-isometric to 
products of finitely many trees 
not only contains 
products of virtually free groups 
but also 
the finitely presented 
simple groups 
constructed by 
Burger and Mozes; 
see~\cite{fp-simpleGs+prod(trees)}.  
% 

%A uniform bound 
%on prime factors 
%of scale function values 
%on envelopes
%of $\Gamma$ 
%is  
%a weak finiteness result 
%for envelopes 
%of $\Gamma$. 

\section{Conventions and outline of the paper}

%We begin with 
%some general remarks. 
%
In what follows, 
$\Gamma$ will denote 
a group 
which is virtually free 
of finite rank 
at least $2$,  
while 
$F$ will denote 
a free group 
of finite rank 
at least $2$. 
We will 
first 
treat 
the case of free groups 
%of finite rank 
%at least $2$ 
to obtain Theorem~\ref{thm:volume_bound(freeG-envelope)}, 
which gives more precise information 
in this special case.  
Then 
we will deduce 
the announced result 
for virtually free groups, 
Corollary~\ref{cor:bound(div(scales))_virt-freeGs}, 
as a corollary. 

To ease discussion 
of the circle of ideas 
which form the topic of this paper, 
we adapt 
the following terminology.  
Let $\Lambda$ and $G$  
be a locally compact groups.  
In our applications, 
$\Lambda$ will usually be 
finitely generated 
and discrete. 
An injective homomorphism 
$\varphi\colon \Lambda \to G$ 
such that 
$\varphi(\Lambda)$ is 
a closed, cocompact subgroup 
in $G$ 
will be called an \emph{envelope of $\Lambda$}. 
The group $G$ 
will also be called 
an envelope of $\Lambda$ 
by abuse of language; 
any envelope 
of a compactly generated group 
is compactly generated. 

In Section~\ref{sec:quasi-action->action} 
we will use 
a quasi-isometric rigidity result 
by 
Lee Mosher, Michah Sageev, and Kevin Whyte 
to see that 
any envelope of $\Gamma$ 
acts cocompactly on a locally finite tree. 
From this  
we deduce in Section~\ref{sec:tree-envelopes} 
that,  
in terms of scale values, 
$\Gamma$ has 
a larger envelope 
which is an automorphism group 
of a locally finite tree. 
Specializing 
to free groups 
in Section~\ref{sec:Schottky-bounds} 
and using 
Lubotzky's results 
on Schottky groups 
of automorphisms 
of trees, 
we bound 
the geometry of 
the underlying trees. 
The required results 
are then 
immediate consequences. 

In this note 
we adhere to 
the following conventions: 
$0$ is a natural number.  
Graphs 
may have loops and multiple edges, 
but edges 
will not be given 
an orientation, 
except 
on occasion of
applications  
of the Bass-Serre theory 
of groups acting on trees.  
The reader 
may consult 
\cite{trees}, section~2.1 
for a formal setup 
of the terminology 
and notation 
for graphs.    
An isometry of a tree 
will be called 
elliptic or hyperbolic 
according to 
whether it admits a fixed point 
(which may be a geometric edge) 
%in the tree 
or not.  
Conjugation by a group element $g$ 
is understood to be 
the map $x\mapsto gxg^{-1}$. 
The relations  $\subset$, $\vartriangleleft$ {\it  etc.\/} always imply
strict inclusion. 
Any automorphism of a topological group will be assumed to be a homeomorphism.

\section{Envelopes of virtually free groups of finite rank act on bushy trees}
\label{sec:quasi-action->action}

Any envelope $G$ 
of a finitely generated group $\Lambda$ 
quasi-acts 
on each of the Cayley graphs 
of $\Lambda$. 
If $\Gamma$ is 
virtually free 
of finite rank 
at least $2$, 
more can be said, 
thanks to 
the following theorem 
by Lee Mosher, Michah Sageev, and Kevin Whyte.  
The theorem 
shows, 
that 
any envelope $G$ 
of $\Gamma$ 
is 
a compact extension of 
a cocompact subgroup of 
the automorphism group 
of a locally finite bushy tree $T$ 
(a tree $T$ 
is \emph{bushy} 
if 
each point of $T$ 
is a uniformly bounded distance 
from a vertex
having at least $3$ 
unbounded complementary components).  
This will enable us 
to reduce to 
the case where 
the envelope $G$ 
is 
the group $\Aut T_.$, 
for such a tree $T$ 
in the next section. 

%\subsection{}

\begin{theorem}[{\cite[Theorem~9]{quasi-Act>T-bound_val}}]% (Virtually free, cocompact lattices)]
\label{thm:virtF,cocp latt}
Let $G$ be a locally compact topological group 
which contains 
a cocompact lattice 
which is 
virtually free of finite rank 
at least $2$. 
Then there exists a cocompact action of $G$ 
on a bushy tree $T$ 
of bounded valence,
inducing 
a continuous, closed homomorphism $\tau\colon G\to \Aut T_.$  
with compact kernel
and cocompact image.
\end{theorem} 
%\pagebreak 

\begin{remark}\label{rem:virtF,cocp latt}
~\\[-4ex]
\begin{enumerate}
\item 
\label{rem:virtF,cocp latt(1)}
The statement 
of Theorem~9 
on page~125 
in~\cite{quasi-Act>T-bound_val} 
uses the word `proper'
in place of `closed'. 
The proof 
of Theorem~9
on page~161 
in~\cite{quasi-Act>T-bound_val} 
makes it clear, 
that 
the intended meaning 
of proper 
in this context 
is closed. 
\item
\label{rem:virtF,cocp latt(2)}
In Theorem~\ref{thm:virtF,cocp latt} 
we can assume that 
the action of $G$ on $T$ via $\tau$ 
is minimal, 
replacing $T$ with 
the minimal $G$-invariant subtree 
if necessary. 
Indeed, 
the $G$-action 
on the tree 
constructed 
in the course of 
the proof of 
Theorem~\ref{thm:virtF,cocp latt} 
in \cite{quasi-Act>T-bound_val}
is already minimal. 
As we are only interested 
in groups containing lattices of positive rank, 
such a minimal tree 
will have no vertex of degree $1$. 
\item 
\label{rem:virtF,cocp latt(3)}
By Theorem~16 
in~\cite{analog(CayleyGphs(topGs))}, 
totally disconnected envelopes  
of a group 
which is 
virtually free 
of finite rank 
can be characterized as 
those groups 
whose rough Cayley graphs 
are quasi-isometric 
to a tree. 
That tree 
is bushy 
if and only if 
the rank 
of free subgroups 
of the cocompact discrete subgroup 
are at least $2$.  
\end{enumerate}
\end{remark}

\begin{corollary}\label{cor:tdlcG-envelope(freeG)}
Let $G$ be 
an envelope 
of a virtually free group 
of finite rank 
at least $2$. 
Then the connected component of the identity 
of $G$ is compact. 
\end{corollary}
\pf{}
By Theorem~\ref{thm:virtF,cocp latt}, 
$G$ has 
a continuous
%, closed 
homomorphism $\tau$ 
into the automorphism group 
of a locally finite tree $T$ 
with compact kernel. 
Since 
$\Aut T_.$ is totally disconnected,  
the connected component $G_0$ of the identity in $G$ 
is contained in 
the kernel of $\tau$. 
We conclude that 
$G_0$ is relatively compact, 
being contained in $\ker(\tau)$.  
Since $G_0$ is closed as well, 
we conclude that 
$G_0$ is compact 
as claimed. 
\qed 

\section{Scale function on envelopes which are automorphism groups of bushy trees}
\label{sec:tree-envelopes}
%Let $F$ be 
%a free group 
%of finite rank 
%at least $2$. 
%
The next three results 
will show,  
that values of the 
scale function 
on totally disconnected 
envelopes of 
the free group $F$ 
can be bounded 
in terms of 
values of the scale function 
on envelopes 
of $F$ 
which are 
automorphism groups  
of some locally finite, bushy tree $T$. 

First, 
we show that 
whenever $G$ is 
an envelope of $F$, 
then 
the codomain of the map $\tau$, 
introduced in Theorem~\ref{thm:virtF,cocp latt}, 
is also an envelope of $F$. 

\begin{proposition}\label{prop:larger-envelop}
Let $\varphi\colon F\to G$ be 
an envelope of 
a free group $F$ 
of finite rank 
at least $2$. 
Let $\tau$ be 
the homomorphism 
from $G$ 
into the automorphism group 
of the tree $T$ 
provided by Theorem~\ref{thm:virtF,cocp latt} 
and 
let $\theta\colon F \to \Aut T_.$ 
be the composite of $\varphi$ and $\tau$. 
Then 
$\theta$ is injective 
and 
$\theta(F )$ is 
a cocompact lattice in $\Aut T_.$. 
Hence $\Aut T_.$ 
is an envelope of $F $. 
\end{proposition}
\pf{}
The group $\theta(F )$ 
is a discrete subgroup of $\Aut T_.$, 
because the kernel of $\tau$ is compact.  
The kernel of $\theta$ is 
a compact subgroup 
of the discrete group $F $, 
hence is finite.  
Since we are assuming that 
$F $ is a free group, 
it is torsion free 
and we conclude that 
the kernel of $\theta$ 
is trivial 
and $\theta$ is injective.  

Since $F $ is cocompact in $G$ 
and $\tau(G)$ is cocompact in $\Aut T_.$ 
we conclude that 
$\theta(F )$ is cocompact in $\Aut T_.$. 
We have verified 
all parts of our claim. 
\qed 

The next lemma 
shows that 
the scale
of an element 
does not change 
if we apply 
a homomorphism, 
which is 
a perfect map
(which, 
for a homomorphism 
of groups,  
means 
a continuous, open, surjective homomorphism 
with compact kernel). 

\begin{lemma}\label{lem:proper_maps&scale}
Let $\pi\colon G\to \widehat{G}$ 
be a continuous, open, surjective homomorphism 
with compact kernel 
between totally disconnected, locally compact groups. 
Let $\alpha$ be 
an automorphism 
of $G$ 
preserving $\ker(\pi)$  
and 
$\widehat{O}$ a subgroup of $\widehat{G}$ 
tidy for 
the automorphism $\widehat{\alpha}$ 
induced by $\alpha$ on $\widehat{G}$. 
Then 
the group 
$\pi^{-1}(\widehat{O})$ 
is tidy for $\alpha$ 
and 
$s_G(\alpha)=s_{\widehat{G}}(\widehat{\alpha})$. 
\end{lemma}
\pf{}
Put $O:=\pi^{-1}(\widehat{O})$. 
Since $\pi$ is continuous 
and 
the kernel of $\pi$ is compact, 
the group $O$ 
is compact and open. 
The conclusion follows from 
the definition of the scale of $\alpha$ 
as a minimum, 
the equation 
\[
|\alpha(O)\colon \alpha(O)\cap O|=
|\widehat{\alpha}(\widehat{O})\colon \widehat{\alpha}(\widehat{O})\cap \widehat{O}|=
s_{\widehat{G}}(\widehat{\alpha}) 
\]
and Proposition~4.7 in~\cite{furtherP(s(tdG))}. 
\qed 

The larger envelope 
of $F$ 
obtained from Proposition~\ref{prop:larger-envelop} 
does not have 
smaller scale values, 
as shown by 
the next result. 

\begin{corollary}\label{cor:upper_bound-Aut(T)}
Let $G$ be 
a totally disconnected, locally compact group, 
$T$ a locally finite, bushy tree  
and 
$\tau\colon G\to \Aut T_.$ 
a continuous, closed homomorphism 
with compact kernel. 
Then 
\[
s_G(x)\leq s_{\Aut T_.}(\tau(x))\quad \text{for all }x\in G \,.
\]
\end{corollary}
\pf{}
Let $x$ be an element of $G$. 
Put $\widehat{G}:=\tau(G)$,  
let $\pi\colon G\to \widehat{G}$ 
be the map 
induced by $\tau$ 
and 
let $\alpha$ be conjugation by $x$. 
Then 
$\pi$ and $\alpha$
satisfy 
the conditions 
on the maps 
with the same names 
in Lemma~\ref{lem:proper_maps&scale} 
and 
$\widehat{\alpha}$ 
is conjugation by $\pi(x)$. 
We conclude that 
$s_G(x)= s_{\widehat{G}}(\pi(x))$. 

Furthermore, 
renaming $\alpha$ as 
conjugation by $\tau(x)$ 
and 
applying 
Proposition~4.3 in \cite{furtherP(s(tdG))}
with $H:=\widehat{G}$, 
we deduce that 
$s_{\widehat{G}}(\pi(x))\le s_{\Aut T_.}(\tau(x))$. 

Combining 
these two relations 
we obtain 
$s_G(x)\le s_{\Aut T_.}(\tau(x))$. 
Since the choice of $x$ was arbitrary, 
the claim follows. 
\qed 

The scale function 
of a closed subgroup, 
$G$,  
of the automorphism group 
of a locally finite tree, 
$T$,  
can be determined geometrically, 
as seen in 
the next result, 
Lemma~\ref{lem:scale(hyp-iso(tree)); prod-formula}.

The geometric description 
of the 
value of 
the scale function 
with respect to $G$ 
at a hyperbolic isometry, 
$h$ say, 
in $G$ 
given in Lemma~\ref{lem:scale(hyp-iso(tree)); prod-formula}  
uses 
ramification indices 
of a subtree $T_{G,\epsilon}$ 
of $T$, 
which 
depends on 
the group $G$ 
and 
the attracting end, 
$\epsilon$, 
of $h$. 
The tree $T_{G,\epsilon}$ 
is defined 
as follows: 
Given an end $\epsilon$ 
of $T$ 
the tree $T_{G,\epsilon}$\label{T_G,epsilon} 
is the union 
of the axes 
of all hyperbolic isometries in $G_\epsilon$.

\begin{lemma}[{\cite[Lemmas~26 and~31]{direction(aut(tdlcG))}}]
\label{lem:scale(hyp-iso(tree)); prod-formula}
Let $G$ be a closed subgroup of the automorphism group of 
a locally finite tree $T$. 
\begin{enumerate}
\item
Let $g$ be an elliptic isometry of $T$ in $G$.   
Then $g$ is topologically periodic 
and $s_G(g)=1=s_G(g^{-1})$. 
\item
Let $h$ be a hyperbolic isometry of $T$ 
in $G$ 
with attracting end $\epsilon$,  
of translation length $n$ say. 
Let  $q_1+1,\ldots,q_n+1$ 
be the ramification indices of 
$n$ consecutive vertices 
on the axis of $h$ 
with respect to the tree $T_{G,\epsilon}$.  
Then $s_G(h)= \prod_{i=1}^n q_i$. 
%
%In particular,  
%for $\epsilon$ in  $\partial_G X$  
%the scale function of $G$ 
%is constant 
%on the cosets of $G_\epsilon$ 
%with respect to $\mathcal{E}G_\epsilon$.  
\end{enumerate}
\end{lemma}

The geometric description 
of the scale function 
given by 
Lemma~\ref{lem:scale(hyp-iso(tree)); prod-formula} 
will be useful 
in the next section.

\section{Upper bounds for scales on envelopes of Schottky lattices} 
\label{sec:Schottky-bounds}

%\cite{nonuT-lattices}

The isomorphic image 
$\theta(F )$ 
of $F$ 
obtained 
in Proposition~\ref{prop:larger-envelop},  
is 
a finitely generated, torsion free, discrete subgroup  
of $\Aut T_.$. 
A finitely generated, torsion free, discrete subgroup 
of the automorphism group of a tree 
will be called 
a \emph{Schottky subgroup of $\Aut T_.$}.  
This terminology 
was coined by 
Alexander Lubotzky 
in section~1 
of the paper~\cite{l.rank1},  
where he describes 
the structure 
of such groups. 

Let $\Gamma$ be a Schottky subgroup of $\Aut T_.$. 
Then $\Gamma$ 
acts freely on $T$ 
and hence 
is a free group;  
%and 
it has finite rank, 
because it is 
finitely generated 
by assumption. 
Bass-Serre theory 
(\cite{trees}, chapter~I or \cite{CovGG}) 
provides us 
with a nice basis of $\Gamma$,  
via its identification with 
the fundamental group 
of the trivial graph of groups 
on the quotient graph $X:=\Gamma \backslash T$.

Any basis  
of a Schottky subgroup, 
$\Gamma$, 
of $\Aut T_.$ 
that is  
obtained from 
Bass-Serre theory 
in the manner explained below 
will be called a \emph{Schottky basis} 
in what follows. 
We recall now 
how these bases 
of $\Gamma$ 
are obtained, 
following the proof 
of Proposition~1.7 in~\cite{l.rank1}. 
This 
construction 
of Schottky bases 
for the group $\Gamma$ 
will involve 
several choices, 
only some of which 
will be of interest 
for us 
later; 
see Lemma~\ref{lem:essential-prop(Schottky_bases)}.  

Choose 
an orientation 
of the edges of $X$. 
Any choice 
of a maximal subtree, 
$Y$ say, 
in $X:=\Gamma\backslash T$ 
defines 
the set of edges $X\smallsetminus Y$, 
belonging to $X$ 
but not to $Y$,   
which freely generate 
the fundamental group 
of the graph $X$. 
The set of edges $X\smallsetminus Y$ 
%so obtained 
corresponds to 
some set of isometries 
of the universal covering tree $T$ 
of $X$. 
This correspondence  
is defined 
in two stages. 

In the first stage 
defining this correspondence, 
choose 
connected subgraphs 
$Y_T\subseteq X_T$ 
in $T$ 
such that 
the canonical projection, 
$p\colon T\to X$,  
is bijective on 
the set of vertices of $Y_T$ 
and 
the set of edges of $X_T$. 
Make this choice 
in such a way 
that the origins of 
all the edges of $X_T\smallsetminus Y_T$ 
belong to $X_T$. 
Such a pair $(Y_T,X_T)$ 
is called 
a `\emph{lifting}' or `\emph{opening}' 
of $X$ 
in the literature.  %\cite[section~5.4]{trees}  
For later use, 
denote 
the inverse image 
of a vertex,
$v$, 
in $Y$ 
under $p$ 
that is 
contained in $Y_T$ 
by $\widetilde{v}$ 
and 
the inverse image 
of an edge,
$e$, 
in $X$ 
under $p$ 
that is 
contained in $X_T$ 
by $\widetilde{e}$.   
Also, 
for an edge, 
$e$, 
denote by 
$o(e)$ and $t(e)$ 
the origin and terminal vertices 
of $e$. 
For the second stage 
defining 
the set of 
isometries 
corresponding to $X\smallsetminus Y$ 
choose, 
for each edge $e$ 
in $X\smallsetminus Y$,  
an element, 
$\gamma_e$, 
of $\Gamma$ 
such that, 
$t(\widetilde{e})=\gamma_e.\widetilde{t(e)}$.  

%  %$d(\gamma_e.\widetilde{o},\widetilde{t})=1$. 
This construction 
is illustrated 
in a simple example 
below, 
where $\Gamma$ 
is the free group 
on~$3$ generators. 
To the left, 
we see 
some possible 
quotient graph $\Gamma\backslash T$. 
The choice 
of a maximal subtree, 
$Y$, 
in $\Gamma\backslash T$ 
is indicated by 
labeling 
the edges 
outside $Y$ 
by the symbols $e_1$, $e_2$ and~$e_3$. 
To the right, 
we see 
how 
the graph $X_T$ 
of an opening 
of the graph $\Gamma\backslash T$ 
determined by 
the choice 
of orientation 
and maximal subtree 
looks like. 
Only the orientation 
of the edges 
$e_1$, $e_2$ and $e_3$ 
matters, 
and 
is indicated by 
arrows 
in both pictures. 
The vertices 
of the graph $X_T$, 
that do not belong to 
the graph $Y_T$, 
are indicated by 
hollow circles. 
All other vertices 
in both pictures 
are indicated by 
filled circles. 
%

% the quotient 
%\fbox{
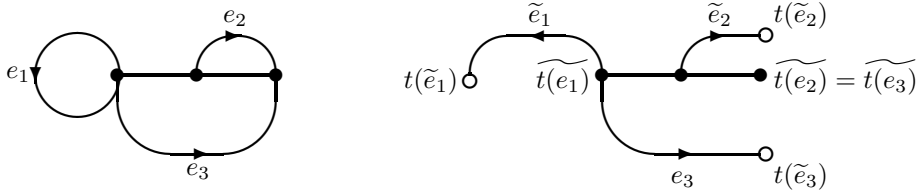
\begin{picture}(160,80)(-10,10)
\thicklines
\put(30,55){\circle{30}}
\put(14,57){\vector(0,-1){5}}
\put(3,54){$e_1$}

\put(45,55){\circle*{5}}
\put(45,55){\line(1,0){30}}
\put(75,55){\circle*{5}}
\put(75,55){\line(1,0){30}}
\put(105,55){\circle*{5}}

\put(75,55){\oval(60,60)[b]}
\put(75,25){\vector(1,0){5}}
\put(71,17){$e_3$}

\put(90,55){\oval(30,30)[t]}
\put(89,70){\vector(1,0){3}}
\put(85,75){$e_2$}
\end{picture}
%}
% the opening 
%\fbox{
\begin{picture}(190,80)(0,15)
\thicklines
\put(0,55){$t(\widetilde{e}_1)$}
\put(25,58){\circle{5}}
\put(50,60){\oval(50,30)[t]}
\put(52,75){\vector(-1,0){5}}
\put(47,80){$\widetilde{e}_1$}

\put(50,55){$\widetilde{t(e_1)}$}
\put(75,60){\circle*{5}}
\put(75,60){\line(1,0){30}}
\put(105,60){\circle*{5}}
\put(105,60){\line(1,0){30}}
\put(135,60){\circle*{5}}
\put(140,55){$\widetilde{t(e_2)}=\widetilde{t(e_3)}$}

\put(105,60){\oval(60,60)[bl]}
\put(105,30){\line(1,0){30}}
\put(105,30){\vector(1,0){5}}
\put(101,19){$e_3$}
%\put(101,22){$e_3$}
%\put(79,37.8){\vector(1,-1){5}}
%\put(70,30.8){$\widetilde{e}_3$}
\put(137,30){\circle{5}}
\put(140,19){$t(\widetilde{e}_3)$}

\put(135,60){\oval(60,30)[tl]}
\put(119,75){\vector(1,0){3}}
\put(137,75){\circle{5}}
\put(115,80){$\widetilde{e}_2$}
\put(140,80){$t(\widetilde{e}_2)$}
\end{picture}
%}
\\

The set $\{\gamma_e\colon e\in X\smallsetminus Y\}$, 
which we obtain 
from all the choices 
we made, 
is a free set 
of generators 
of~$\Gamma$ 
by \cite[\S 5]{trees}, 
which we call 
the Schottky basis of $\Gamma$ 
determined by these choices.

The Schottky basis $\{\gamma_e\colon e\in X\smallsetminus Y\}$ 
consists of 
hyperbolic isometries of $T$. 
This is true 
because 
the set of edges 
of the graph $X_T$  
is a fundamental domain 
for the $\Gamma$-action 
on the edges. 
The latter argument 
works 
for any set $\{\gamma_e\colon e\in X\smallsetminus Y\}$
of elements 
of a group $\Gamma$ 
constructed 
in the way 
described above, 
as long as $\Gamma$ 
acts without inversion 
of edges.  
In the case at hand, 
where $\Gamma$ 
is a free group, 
it may 
alternatively 
be seen, 
by using that 
$\Gamma$ must act 
freely on $X$. 

Below, 
we will need 
the following information 
on the axes 
(and translation lengths) 
of elements 
in the Schottky basis $\{\gamma_e\colon e\in X\smallsetminus Y\}$. 
The axis 
of the element $\gamma_e$ 
passes through 
the vertices $\widetilde{t(e)}$ and $t(\widetilde{e})=\gamma_e.\widetilde{t(e)}$, 
whose distance 
is therefore 
the translation length 
of $\gamma_e$. 
This may be seen 
by applying 
Lemma~1.2 
in~\cite{l.rank1} 
with 
$\gamma$ equal to $\gamma_e$, 
$x$ equal to $\widetilde{t(e)}$ 
and 
$y$ equal to the first vertex 
on the shortest path 
joining $\widetilde{t(e)}$ to $t(\widetilde{e})$, 
(which, by the way, 
is contained in $X_T$).  
The axis of $\gamma_e$ 
also passes through $o(\widetilde{e})$, 
because 
by our choice 
of a lift for $e$, 
the vertex $t(\widetilde{e})$  
is incident to 
only one edge of $X_T$.  
%\cite{E(lattices<Kac-MoodyG(finiteF)) 

This finishes 
the description 
of the construction of 
the set of 
Schottky bases for $\Gamma$, 
save for 
two final remark 
for readers 
familiar with~\cite{l.rank1}. 
First, 
in~\cite{l.rank1}, 
a Schottky basis 
is supposed to 
satisfy 
an apparently stronger condition 
than the one 
we have verified here, 
see Definition~1.4 in~\cite{l.rank1}. 
But, 
as stated in 
the proof 
of Propostion~1.7 in~\cite{l.rank1}, 
the required labeling 
of the axes of 
the elements 
in the basis $\{\gamma_e\colon e\in X\smallsetminus Y\}$ 
is obtained 
if we choose, 
for every $e\in X\smallsetminus Y$, 
the labeling 
of the axis of $\gamma_e$ 
determined by 
labeling $\widetilde{t(e)}$ 
with the symbol $x_1$. 
Second, 
the construction above 
gives 
all Schottky bases 
for $\Gamma$ 
in the sense 
defined in 
Definition~1.4 in~\cite{l.rank1}. 
To see this, 
note that 
the fundamental domain~$F$ 
from Proposition~1.6 in~\cite{l.rank1} 
defines an opening 
for $\Gamma\backslash T$, 
such that 
the given Schottky basis $\{\gamma_1,\ldots,\gamma_l\}$ 
is one of the possible Schottky bases 
obtainable from 
that opening 
in the sense 
and the way 
described 
above. 

The next result 
is obvious from 
the discussion 
of the construction 
of a Schottky basis 
above. 

\begin{lemma}\label{lem:essential-prop(Schottky_bases)}
Let $T$ be a tree. 
Suppose that 
$\{\gamma_1,\ldots,\gamma_l\}$ 
and 
$\{\gamma_1',\ldots,\gamma_l'\}$ 
are 
two Schottky bases of 
the same Schottky subgroup, 
$\Gamma$,  
in $\Aut T_.$, 
obtainable from 
the same choice 
of maximal subtree, 
$Y$ say, 
in $\Gamma\backslash T$. 
Let $\{e_1,\ldots, e_l\}$ 
be the set of edges 
of $\Gamma\backslash T$ 
outside $Y$. 
Then there are 
permutations, 
$\pi$ and $\pi'$,  
of the set 
$\{1,\ldots,l\}$ 
such that 
for each integer $i$ 
with $1\le i\le l$ 
the restriction of 
the canonical projection $p\colon T\to \Gamma\backslash T$  
induces bijections 
from  
$(1)$ onto~$(0)$ 
and 
from 
$(1')$ onto~$(0)$  
where: 
\begin{itemize}
\item[$(0)$] 
is 
the set 
of vertices 
on the shortest path 
in 
%the tree 
$Y$ 
connecting 
the two vertices 
of the edge $e_i$; 
\item[$(1)$] 
is 
the set 
of vertices 
of a fundamental domain 
for the vertices 
on the axis 
of $\gamma_{\pi(i)}$ 
modulo $\langle\gamma_{\pi(i)}\rangle$; 
\item[$(1')$] 
is 
the set 
of vertices 
of a fundamental domain 
for the vertices 
on the axis 
of $\gamma'_{\pi'(i)}$ 
modulo $\langle\gamma'_{\pi'(i)}\rangle$. 
\end{itemize}
\end{lemma}
%\pf{}
%\qed 

Let $T$ be a tree 
and 
$\Gamma$ be 
a Schottky subgroup of $\Aut T_.$.  
In what follows, 
we will only be interested in 
those properties 
of a Schottky basis,  
that 
are already determined by 
the quotient graph $\Gamma\backslash T$ 
and 
the set of 
collections  
of vertices described in~(0), 
indexed by 
the set  $\{e_1,\ldots,e_l\}$ of edges 
outside a maximal subtree, 
$Y$, 
of $\Gamma\backslash T$. 
Therefore, 
we will allow ourselves 
to speak, 
on occasion, 
somewhat inaccurately,  
of 
`the Schottky basis determined by the maximal subtree $Y$'.

We now 
return to 
the Schottky subgroup $\theta(F )$ 
of $\Aut T_.$ 
obtained 
in Proposition~\ref{prop:larger-envelop}.  
The group $\theta(F )$ 
is also a (cocompact) lattice 
in $\Aut T_.$. 
It follows that 
the quotient graph 
$X:=\theta(F )\backslash T$ 
is finite. 
Using 
part~\ref{rem:virtF,cocp latt(2)} 
of Remark~\ref{rem:virtF,cocp latt}, 
we may assume that 
$T$ has no vertex of degree $1$. 
This implies that 
$X$ has no vertex of degree $1$ either. 

We will call a Schottky subgroup of $\Aut T_.$ 
that is also a lattice 
a \emph{Schottky lattice}.
We now derive our main lemma on Schottky lattices. 

\begin{lemma}\label{lem:ramification_bound(SchottkyGs)}
Any tree $T$ % necessarily minimal  
with all vertices 
having degree at least $3$, 
which admits 
a Schottky lattice on $n$ generators, 
belongs to a finite list of 
universal covering trees  
of finite graphs. 
In particular 
there are integers $K_n$ and $l_n$, 
such that 
\begin{enumerate}
\item
the ramification indices of $T$ 
are at most $K_n$; 
\item 
every member of a Schottky basis 
of a Schottky lattice on $T$ 
has translation length at most $l_n$. 
\end{enumerate}
\end{lemma}
\pf{}
It suffices to show, 
that there are 
only finitely many graphs 
with fundamental group 
free of rank~$n$ 
with all vertices of degree at least $3$. 
We will show this,  
by deriving 
upper bounds 
for the number 
of vertices 
and edges 
of such a graph. 

Let $X$ be such a graph. 
It has 
$n$ edges outside a maximal subtree. 
Denote by  
$e$ the total number of edges 
and by 
$v$ the number of vertices 
of $X$.  
Since 
the number of vertices 
of a tree 
exceeds 
the number of its edges 
by one,  
we have $n=e-v+1$. 

The sum 
of the degrees 
of vertices 
in $X$ 
is at most $2e$, 
because 
at most $2$ 
different vertices 
belong to any edge.  
On the other hand 
that sum 
is at least $3v$ 
by our assumption 
on the vertex degrees. 
It follows that 
$e\geq 3/2\,v$.  
Substituting 
this inequality  
into the equation 
stated in the last paragraph 
we obtain 
the upper bound $v\le 2(n-1)$
for $v$ 
in terms of $n$. 

The equation $n=e-v+1$ 
implies that 
$e=n+v-1$, 
which is seen to be 
bounded above 
by $3(n-1)$ 
using the upper bound $v\le 2(n-1)$,  
derived above.  
This 
completes the proof. 
\qed

Finally, 
we derive our main result. 

\begin{theorem}
\label{thm:volume_bound(freeG-envelope)}
Let $F$ be a free group 
of finite rank 
at least $2$. 
There is a constant $V(F )$ 
such that 
for all envelopes 
$\varphi\colon F \to G$ 
of $F $ 
into a totally disconnected group, 
there is a free generating set $S_\varphi$ 
in $F $ 
such that 
\[
\prod_{t\in S_\varphi} s_G(\varphi(t))\leq V(F )\,.
\]
Furthermore, 
the set of primes 
dividing one of the numbers 
\[
s_G(x);\quad \text{$G$ a $F $-envelope, $x\in G$} %F $}
\]
is finite. 
\end{theorem}
\pf{}
Proposition~\ref{prop:larger-envelop} 
and 
Corollary~\ref{cor:upper_bound-Aut(T)} 
show that, 
for both claims, 
we may assume that 
$G$ runs through 
automorphism groups 
of bounded valence, bushy trees, 
which are envelopes of $F$. 
As pointed out 
in Remark~\ref{rem:virtF,cocp latt}, 
we may assume that 
the corresponding 
tree actions 
are minimal, 
hence that 
none of these trees has 
a vertex of degree $1$. 
Since 
$F$ is torsion free, 
it acts 
on the trees 
in an orientation preserving way. 
Hence,  
we may also assume that 
none of the trees 
has vertices  
of degree $2$. 

As explained 
at the beginning of this section, 
this implies that 
the image of $F$ under $\varphi$ 
is a Schottky lattice. 
We choose $S_\varphi$ 
to be a Schottky basis 
of $\varphi(F)$. 
It follows then 
from 
Lemma~\ref{lem:scale(hyp-iso(tree)); prod-formula} 
and 
Lemma~\ref{lem:ramification_bound(SchottkyGs)} 
that 
we have 
$V(F )\leq (K_n-1)^{nl_n}$ 
and 
that no prime divisor of $s_{\Aut T_.}(x)$ 
can exceed $K_n$. 
This shows our claim. 
\qed 

The second property 
in the above theorem 
can also be obtained 
for groups 
which are just virtually free 
of finite rank at least $2$. 

\begin{corollary}\label{cor:bound(div(scales))_virt-freeGs}
Let $\Gamma$ be 
virtually a free group of finite rank 
at least $2$.  
Then 
the set of primes 
dividing one of the numbers 
\[
s_G(x);\quad \text{$G$ a totally disconnected $\Gamma $-envelope, $x\in G$} %F $}
\]
is finite. 
\end{corollary}
\pf{}
Let $F$ be 
a free subgroup 
of finite index 
in $\Gamma$ 
and 
let $\iota\colon F\hookrightarrow \Gamma$ 
be the inclusion of $F$ in $\Gamma$. 
If 
$\varphi\colon \Gamma\to G$ 
is a $\Gamma$-envelope, 
then 
the composite $\varphi\circ\iota\colon F\to G$ 
is an $F$-envelope. 
Therefore 
the claim 
follows immediately 
from the second statement 
of Theorem~\ref{thm:volume_bound(freeG-envelope)}. 
\qed

We now 
make 
terminology 
and 
notation 
available, 
to describe 
the phenomenon 
established 
for free groups 
of finite rank 
at least~$2$ 
in 
the first part of 
Theorem~\ref{thm:volume_bound(freeG-envelope)}. 
This 
terminology 
and 
notation 
will be used 
in the next section. 

\begin{definition}\label{def:scale-volumes}
Let $\Lambda$ be a finitely generated group 
and 
let $\varphi\colon \Lambda\to G$ 
be an envelope of $\Lambda$ 
with 
totally disconnected, locally compact codomain $G$. 
\begin{enumerate}
\item % scale-volume of a generating set 
If $S$ is 
a finite set 
of generators 
for $\Lambda$,  
the number 
$\svol_\varphi(S):=\prod_{t\in S} s_G(\varphi(t))$ 
will be called 
\emph{the scale volume of $S$ with respect to $\varphi$}. 
\item % scale-volume of a embedding 
The number 
$\svol_\varphi\langle\Lambda\rangle:=
\min\{\svol_\varphi(S)\colon \text{$S$ is finite and generates $\Lambda$}\}$ 
will be called 
\emph{the scale volume of $\Lambda$ with respect to $\varphi$}. 
\item % scale-volume of a group
The number 
$\Svol(\Lambda):=
\sup\{\svol_\varphi\langle\Lambda\rangle\colon 
\text{$\varphi$ is an envelope of $\Lambda$ with totally disconnected codomain}\}$ 
(which may be infinite) 
will be called 
\emph{the scale volume of $\Lambda$}.
\end{enumerate}
\end{definition}

Using 
the terminology 
just introduced, 
the first part 
of Theorem~\ref{thm:volume_bound(freeG-envelope)} 
may be 
restated as 
`the scale volume 
of a free group 
of finite rank 
at least $2$ 
is finite'. 

\section{Explicit bounds on scales}\label{sec:explicit bounds}

In this section 
we will address 
the problem 
of getting 
quantitative versions 
of 
Theorem~\ref{thm:volume_bound(freeG-envelope)} 
and 
Corollary~\ref{cor:bound(div(scales))_virt-freeGs}. 
We start with 
explicit bounds 
on possible values 
of the scale function 
on a free group 
of finite rank 
at least $2$.  

\begin{proposition}
\label{prop:explicit_bounds}
Let $n$ be an integer 
which is at least $2$. 
Let $\mathcal{B}(n)$ 
be the class 
of all graphs 
with all vertices of degree at least $3$ 
and fundamental group 
free of rank $n$. 
Then 
\begin{enumerate}
\item 
The maximal vertex degree 
of graphs in $\mathcal{B}(n)$ 
is $2n$. 
This maximal degree is achieved 
for the unique graph, 
$B_{2n}$, 
in $\mathcal{B}(n)$ 
with $1$ vertex 
and 
is not achieved 
for any other graph in $\mathcal{B}(n)$. 
The fundamental group, 
$F$, 
of the trivial graph of groups over $B_{2n}$ 
admits, 
via its action 
on the universal covering tree, 
$T_{2n}$,  
of $B_{2n}$, 
which is a homogenous tree 
of degree $2n$, 
a discrete, cocompact embedding 
into $\Aut {T_{2n}}_.$. 
All elements 
of the unique Schottky basis 
of $F$ 
obtained using all edges of $B_{2n}$ 
have translation length~$1$ 
and 
the scale function of $\Aut {T_{2n}}_.$ 
assumes the value $2n-1$ 
on all of these elements.  
\item 
For any odd integer $s$ 
with $3\le s\le 2n-1$, 
there is a graph, 
$B(s)$,  
in $\mathcal{B}(n)$ 
such that 
the scale function 
on the automorphism group 
of the covering tree, 
$T(B(s))$, 
of $B(s)$ 
assumes the value~$s$ 
on 
some hyperbolic isometry 
of translation length~$1$  
which 
belongs to 
every Schottky basis 
of the fundamental group, 
$F$,  
of the trivial graph of groups 
over $B(s)$. 
The group $F$ 
admits,  
via its action 
on $T(B(s))$,  
a discrete, cocompact embedding 
into $\Aut{T(B(s))}_.$.  
\item 
Let $n$ be an integer 
which is at least $2$. 
The maximal translation length 
of a member of 
a Schottky basis 
for the image 
of the free group 
on $n$ generators 
in a discrete, cocompact embedding 
in the automorphism group 
of a locally finite tree 
all of whose vertices 
have degree 
at least $3$
is $2(n-1)$. 

%
%This maximum is attained. 
More precisely, 
there is a graph, 
$B^{2(n-1)}$,  
in $\mathcal{B}(n)$ 
such that 
every Schottky basis 
of the fundamental group, 
$F$, 
of the trivial graph of groups 
over $B^{2(n-1)}$ 
has an element 
whose translation length is $2(n-1)$ 
and whose scale 
with respect to 
the automorphism group 
of the covering tree, 
$T(B^{2(n-1)})$, 
of $B^{2(n-1)}$ 
is $2^{2(n-1)}$. 
The group~$F$ 
admits,  
via its action 
on $T(B^{2(n-1)})$,  
a discrete, cocompact embedding 
into $\Aut{T(B^{2(n-1)})}_.$.   
%\item 
\end{enumerate}
\end{proposition}
\pf{}
\noindent\textbf{Proof of 1:\quad}
In proving 
the first claim, 
we start by 
confirming  
the statements 
made 
in the first two sentences 
thereof. 
Thanks to 
the lower bound 
on vertex degrees 
for graphs is $\mathcal{B}(n)$, 
the following holds 
for any graph, 
$B$ say, 
in $\mathcal{B}(n)$. 
Contracting an edge 
of a maximal subtree 
in $B$ 
to its initial vertex, 
$o$ say,  
creates another graph 
in $\mathcal{B}(n)$, 
whose vertex degree 
at $o$
is strictly larger than 
the vertex degree 
of $o$ in $B$. 
Therefore, 
the maximal vertex degree 
of a graph in $\mathcal{B}(n)$ 
is achieved 
for a graph 
with just one vertex. 
Since 
there is 
a unique graph, 
$B_{2n}$, 
in $\mathcal{B}(n)$ 
with one vertex, 
we have seen that 
the first two sentences 
are true. 

A moment's thought shows, 
that 
all the remaining statements 
made in the first claim 
save the value 
taken by the scale function 
of $\Aut{T_{2n}}_.$ 
on the elements 
of the Schottky basis 
of $F$ 
follow immediately from 
the statements just proved. 

To see that 
this remaining statement 
is true also, 
it suffices 
by Lemma~\ref{lem:scale(hyp-iso(tree)); prod-formula} 
to show 
that 
for any end $\epsilon$ 
of $T_{2n}$ 
the tree ${(T_{2n}})_{\Aut{T_{2n}}_.,\epsilon}$
(introduced 
in the paragraph 
preceding Lemma~\ref{lem:scale(hyp-iso(tree)); prod-formula},  
page~\pageref{T_G,epsilon})   
is the whole tree $T_{2n}$. 
Since 
for any two distinct ends, 
$\epsilon$ and $\epsilon'$ say,  
in a homogenous tree 
such as $T_{2n}$ 
there is a hyperbolic isometry 
(of translation length $1$) 
whose axis is 
the line joining 
$\epsilon$ and $\epsilon'$, 
we clearly have 
 ${(T_{2n}})_{\Aut{T_{2n}}_.,\epsilon}=T_{2n}$, 
for any end $\epsilon$, 
and we finished 
the proof 
of our first claim.

\vspace{1ex}\noindent\textbf{Proof of 2:\quad}
To obtain 
the statements  
of the second claim, 
we first introduce 
the graph $B(s)$ 
for a given odd integer $s$ 
with $3\le s\le 2n-1$. 
If $s$ is equal to $2n-1$, 
then, 
by the first claim, 
which we already proved, 
the graph $B_{2n}$ 
has 
all the properties 
the graph $B(2n-1)$ 
is required to satisfy.  
Therefore, 
we may assume 
in what follows, 
that $s$ is at most $2n-3$. 

If $s\le 2n-3$, 
the graph $B(s)$ 
has $2$ vertices, 
$v_0$ and $v_1$. 
The vertices $v_0$ and $v_1$ 
are joined by 
one, respectively two edges 
depending on whether 
$s$ is different from or equal to $n-1$.   
No matter what 
the value of $s$ is, 
attach $(s+1)/2\ge 2$ loops 
to the vertex $v_0$. 
If $s\neq n-1$, 
attach $n-(s+1)/2\ge n-(2n-2)/2=1$ loops 
to the vertex $v_1$. 
If $s=n-1$, 
attach $n-n/2-1=n/2-1\ge 1$ loops 
to the vertex $v_1$. 

For each odd integer $s$ 
with $3\le s\le 2n-3$ 
the graph $B(s)$ 
belongs to $\mathcal{B}(n)$  
by construction. 
The degrees $d_0$ and $d_1$ 
of the vertices $v_0$ and $v_1$ 
of $B(s)$ 
are different.  

We claim that 
the scale function 
with respect to 
%the automorphism group 
%of the covering tree $T(B(s))$ 
$\Aut T(B(s))_.$ 
assumes the value $s$ 
on each of the elements of $F$ 
obtained from 
the $(s+1)/2$ loops 
based at $v_0$. 
To see this, 
let $h$ be one of 
the elements of $F$ 
obtained from loops at $v_0$ 
and 
let $\epsilon$ 
be the attracting end 
of $h$. 
We will determine 
the tree $T(B(s))_{\Aut T(B(s))_.,\epsilon}$. 

Note that 
all vertices 
on the axis of $h$ 
are translates of 
a lift of $v_0$ 
and hence have 
degree $d_0$. 
The vertices 
on the axis 
of every hyperbolic isometry 
of $T(B(s))$ 
that has $\epsilon$ 
as attracting end 
must have degree $d_0$ 
also. 
Let $k$ be 
such a hyperbolic isometry. 
By the above observation 
on the degrees 
of vertices 
on the axis of $k$, 
this axis 
maps to 
a closed edge path 
at $v_0$ 
in $B(s)$ 
all of whose edges 
belong to 
the subgraph, $B(v_0)$, 
of $B(s)$ 
consisting of 
the vertex $v_0$ 
and 
all the loops based at $v_0$. 
Conversely, 
every closed edge path 
at $v_0$ 
can be obtained 
as an image of 
a hyperbolic isometry 
of the tree $T(B(s))$ 
fixing the end $\epsilon$. 
We conclude that 
the tree $T(B(s))_{\Aut T(B(s))_.,\epsilon}$ 
equals  
the lift 
of $B(v_0)$ 
in $T(B(s))$ 
that contains $\epsilon$.
In particular, 
the tree $T(B(s))_{\Aut T(B(s))_.,\epsilon}$ 
is regular 
of degree $2(s+1)/2=s+1$. 

Our claim 
on the values of the scale function 
follows from Lemma~\ref{lem:scale(hyp-iso(tree)); prod-formula}.  
Since 
the remaining statements made 
are obvious, 
the proof of 
the second claim 
is complete. 

\vspace{1ex}\noindent\textbf{Proof of 3:\quad}
We first establish 
the upper bound 
on the translation length 
claimed. 
Recall that 
the upper bound 
for the number 
of vertices 
of a graph 
with fundamental group free 
of rank $n$ 
and all vertices 
of degree at least $3$ 
obtained in 
the proof of Lemma~\ref{lem:ramification_bound(SchottkyGs)} 
was $2(n-1)$. 
The longest path 
without backtracking 
inside a maximal subtree 
of a graph 
with fundamental group free 
of rank $n$ 
and all vertices 
of degree at least $3$ 
therefore is $2(n-1)-1$ 
and 
the longest translation length 
an element of a Schottky basis 
in the fundamental group 
over the trivial graph of groups 
over such a graph 
can have 
is $2(n-1)-1+1=2(n-1)$. 
This establishes 
the upper bound 
for the translation length 
claimed. 

To prove 
the remainder 
of the claim, 
it suffices 
to prove 
the existence 
of the graph $B^{2(n-1)}$ 
with the properties stated 
for every integer $n$ 
at least $2$, 
since 
the other statements 
follow 
from this. 
The graph  $B^{2(n-1)}$ 
is constructed as follows. 
Take the $2(n-1)$ numbers 
$0,\ldots ,2n-3$ 
as the set of vertices of $B^{2(n-1)}$. 
The number of these vertices is even. 
Choose, 
once and for all, 
a cyclic order 
on the set 
of vertices 
of $B^{2(n-1)}$. 
Traverse 
the vertices 
of $B^{2(n-1)}$
in that cyclic order 
and connect 
encountered successive vertices 
in an alternating fashion 
by $1$ respectively $2$ edges. 
By construction, 
every vertex 
of the graph $B^{2(n-1)}$ 
has degree $3$. 

Since 
the various Schottky bases 
of the fundamental group 
of the trivial graph of groups 
over $B^{2(n-1)}$ 
are determined 
by the choice 
of a maximal subtree 
in $B^{2(n-1)}$, 
we determine next 
all choices 
of a maximal subtree 
in $B^{2(n-1)}$. 
If $T$ is 
a maximal subtree 
of $B^{2(n-1)}$, 
then 
not every vertex 
in $B^{2(n-1)}$ 
can be connected 
to the previous vertex 
within $T$ 
with respect to 
the chosen cyclic order 
on the set of vertices, 
because then 
$T$ would contain 
a cycle. 
Let $v_0$ be 
a vertex 
of $B^{2(n-1)}$
which is not connected to 
its predecessor, 
$v_{2n-3}$,  
in the chosen cyclic order. 
Then 
the set of edges of $T$ 
forms a path 
in cyclic order 
from $v_0$ to $v_{2n-3}$. 

We therefore see that, 
up to an automorphism  
of $B^{2(n-1)}$,   
the maximal subtrees 
in $B^{2(n-1)}$ 
are of 
at most 
two kinds. 
If $n$ is at least $3$ 
there are two kinds,  
depending on 
whether 
$v_{2n-3}$ is connected to $v_0$ 
by $1$ or by $2$ edges 
inside $B^{2(n-1)}$, 
while 
there is just one kind 
if $n$ is equal to $2$.  
Each of the edges 
connecting $v_{2n-3}$ to $v_0$ 
that do not belong to $T$  
defines an element of 
the Schottky basis 
determined by $T$ 
of translation length $2(n-1)$. 
We have just seen that 
there is 
at least $1$ element 
with such a translation length 
in every Schottky basis 
attached to 
the graph $B^{2(n-1)}$ . 
Since 
the covering tree $T(B^{2(n-1)}$ 
of  $B^{2(n-1)}$ 
is $3$-regular, 
the scale 
of such an element 
with respect to $\Aut T(B^{2(n-1)})_.$ 
is $2^{2(n-1)}$.  
We have shown 
all parts 
of claim~3. 
\qed

The following corollary 
follows from 
%parts~2 and~3 of 
Proposition~\ref{prop:explicit_bounds}. 

\begin{corollary}\label{cor:estimates(scale-val(freeGs))}
Let $n$ be an integer 
which is at least $2$  
and 
let $F_n$ be 
the free group 
of rank $n$. 
Then, 
for every prime number, 
$p$,  
such that $2\le p\le 2n-1$,  
there is 
a discrete, cocompact embedding, 
$\varphi\colon F_n\to G$, 
into a totally disconnected, locally compact group, $G$,  
and an element, 
$h_\varphi$, 
of $F_n$ 
such that 
$p\mid {s_G(\varphi(h_\varphi))}$. 
Furthermore, 
every prime, 
$p$,  
dividing 
the value 
of the scale function 
on an element of $F_n$ 
with respect to any 
envelope 
of $F_n$ 
satisfies $2\le p\le 2n-1$.
\qed
\end{corollary}

From Proposition~\ref{prop:explicit_bounds} 
we also obtain 
%the following 
estimates 
for the scale volume 
of the free group  
of rank 
at least~$2$. 

\begin{corollary}\label{cor:estimates(scale-vol(freeGs))}
Let $n$ be 
an integer 
which is at least $2$, 
and 
let $F_n$ be 
the free group 
of rank $n$. 
Then 
the scale volume, 
$\Svol(F_n)$,
of $F_n$ 
satisfies 
the inequalities 
\[
(2n-1)^n \le \Svol(F_n) \le (2n-1)^{2n(n-1)}\,. 
\]
\end{corollary}
\pf{}
The stated lower bound 
follows from 
part~1
of Proposition~\ref{prop:explicit_bounds}.
The upper bound 
follows from 
parts~1 and~3 
of Proposition~\ref{prop:explicit_bounds}.
\qed

The gap between 
the upper and lower bounds 
in the estimate  
for $\Svol(F_n)$ 
in~Corollary~\ref{cor:estimates(scale-vol(freeGs))} 
grows fast, 
and both bounds are 
very far from 
the truth.  
Some experimentation 
led me 
to believe that, 
asymptotically, 
as $n$ goes to infinity,  
the value of 
$\Svol(F_n)$ 
is %approximately 
$2^{n(2\log_2(n/3)+3)}$. 
However, 
I have been forced 
to revise 
my initial impression 
that 
it is straightforward 
to determine 
the asymptotics 
of the sequence $(\Svol(F_n))_{n\ge 2}$.

%	\bibliographystyle{alpha}
%	\bibliography{%
%	short-Abk,%
%	local,%
%	AAGruppen,%
%	%Abk,%
%	Baeume,% 
%	%GGtheory,% 
%	Gitter,% 
%	NCRaeume,% 
%	%arithmetic,% 
%	buildings,% 
%	diverses,% 
%	%growth,%
%	%measures,% 
%	%qi,% 
%	%quartett,% 
%	%short-Abk,% 
%	%tilings,% 
%	topG%, 
%	%types,% 
%	%unclassified%
%	}
%	\end{document}

\end{document}